\documentclass[12pt]{article}
\usepackage[utf8]{inputenc}
\usepackage[english]{babel}
\usepackage{graphicx}
\usepackage{vmumath}
\usepackage{amsmath}
\usepackage{amsthm}
\usepackage{avo_book}
\usepackage{tikz}
\usepackage{caption}
\usepackage{hyperref}
\usepackage{subcaption}
\usepackage{booktabs}
\usepackage[a4paper,left=15mm,right=15mm,top=30mm,bottom=20mm]{geometry}
\noaftermath
\begin{document}

\bibliographystyle{unsrt}

\title{Enumeration of Labelled and Unlabelled Hamiltonian Cycles\\ in Complete $k$-partite Graphs}

\author{
Evgeniy Krasko \qquad Igor Labutin\qquad   Alexander Omelchenko\\
\small St. Petersburg Academic University\\
\small 8/3 Khlopina Street, St. Petersburg, 194021, Russia\\
\small\tt \{krasko.evgeniy, labutin.igorl, avo.travel\}@gmail.com
}

\maketitle

\begin{abstract}
We enumerate labelled and unlabelled Hamiltonian cycles in complete $n$-partite graphs $K_{d,d,\ldots,d}$ having exactly $d$ vertices in each part (in other words, Tur\'an graphs $T(nd, n))$. We obtain recurrence relations that allow us to find the exact values $b_{n}^{(d)}$ of such cycles for arbitrary $n$ and $d$. 

\end{abstract}

\section{Introduction}

The problem of enumerating Hamiltonian cycles in different classes of graphs is one of the most difficult problems of enumerative combinatorics. Apart from some trivial examples (like Hamiltonian cycles in complete graphs), only a few exact results of this type are known. Due to the inherent complexity of such problems, the efforts of researchers have been largely concentrated on obtaining upper and lower bounds on the numbers of Hamiltonian cycles in different classes of graphs (see \cite{Thomassen_Hamiltonian},\cite{Alon_Hamiltonian},\cite{Hamiltonian_Dirac_graphs},\cite{Hamiltonian_Petersen_graphs},\cite{Hamiltonian_n_Cube}). Even fewer results regarding unlabelled Hamiltonian cycles have been obtained so far.

\begin{figure}[ht]
\centering
	\centering
    	\includegraphics[scale=0.7]{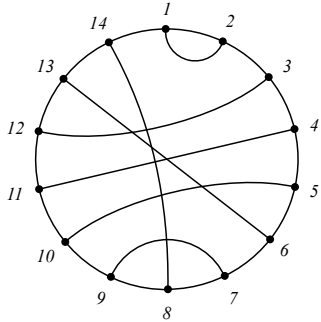}
	\caption{A chord diagram}
\label{fig:chord}
\end{figure}

One exception in this regard is the work \cite{Singmaster} in which the author derived an analytic formula for the numbers $H_n$ of labelled Hamiltonian cycles in $n$-dimensional octahedrons (\url{http://oeis.org/A003436}), that is, in $n$-partite graphs $K_{2,2,\ldots,2}$ having $2n$ vertices. That article also contains a table of the corresponding numbers for unlabelled Hamiltonian cycles, numerically computed for small $n$. 20 years later the numbers $H_n$ appeared once again in the problem of enumerating loopless chord diagrams \cite{Kalashnikov}. A chord diagram consists of $2n$ points on a circle labelled with the numbers $1,2,\ldots,2n$ in a circular order and joined pairwise by chords (figure \ref{fig:chord}). A chord is said to be a loop if it connects two neighboring points (chord $\{1,2\}$ on Figure \ref{fig:chord}). A loopless chord diagram is a diagram without loops. 

\begin{figure}[ht]
\centering
	\begin{subfigure}[b]{0.3\textwidth}
	\centering
    		\includegraphics[scale=0.7]{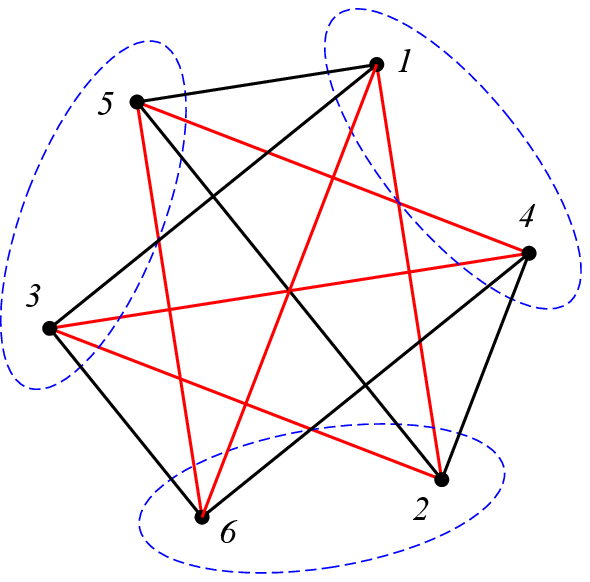}
		\caption{}
	\end{subfigure}	
	\begin{subfigure}[b]{0.3\textwidth}
	\centering
    		\includegraphics[scale=0.7]{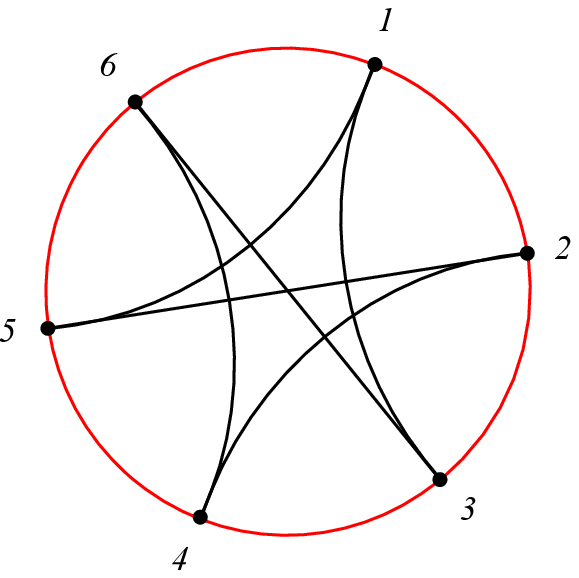}
		\caption{}
	\end{subfigure}
	\begin{subfigure}[b]{0.3\textwidth}
	\centering
    		\includegraphics[scale=0.7]{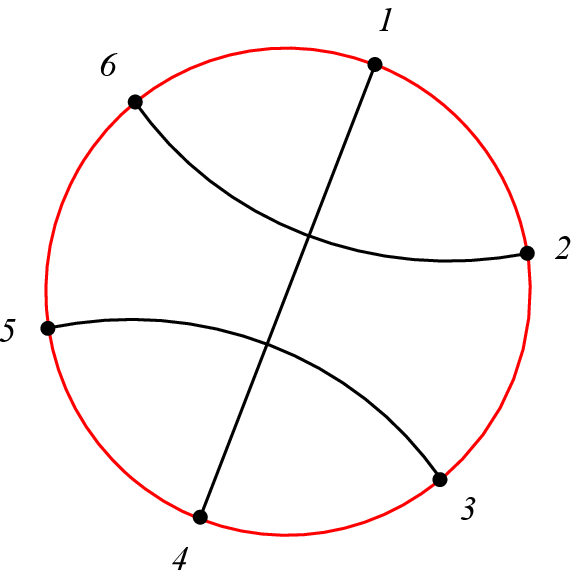}
		\caption{}
	\end{subfigure}
	\caption{Correspondence between Hamiltonian cycles in octahedrons and chord diagrams}
\label{fig:octahedron}
\end{figure}

In the paper \cite{Krasko_Om_chord_diagrams} a bijection between Hamiltonian paths in octahedrons and loopless chord diagrams was noted. Take an $n$-dimensional octahedron with a distinguished Hamiltonian cycle (Figure \ref{fig:octahedron}(a)) and draw it in such a way that this cycle forms a circle on a plane (Figure \ref{fig:octahedron}(b)). Then remove all of its edges that don't belong to the Hamiltonian cycle and add chords between those vertices that weren't connected by an edge before (Figure \ref{fig:octahedron}(c)). The resulting object is a chord diagram which is necessarily loopless: traversing a Hamiltonian cycle in $K_{2,2,\ldots,2}$ we can't visit two vertices of the same part one after another. Clearly, this transformation is invertible.

\begin{figure}[ht]
\centering
	\begin{subfigure}[b]{0.4\textwidth}
	\centering
    		\includegraphics[scale=0.7]{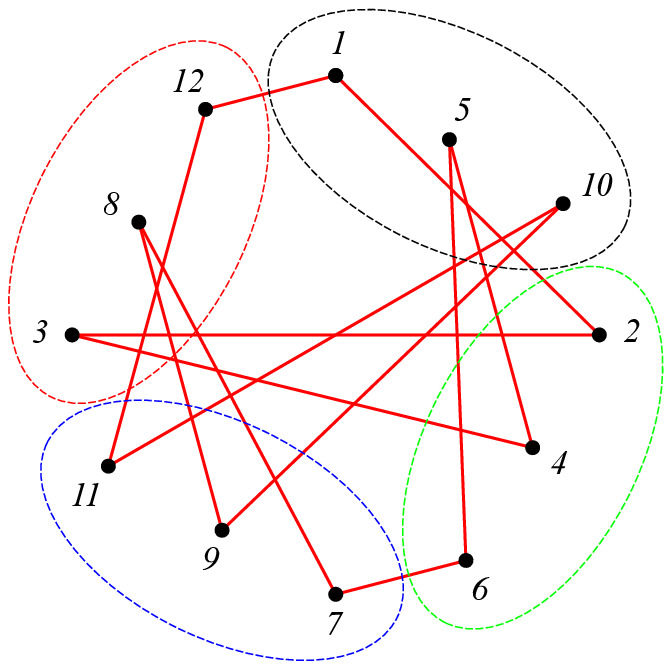}
		\caption{Hamiltonian cycle in $K_{3,3,3,3}$}
	\end{subfigure}
	\begin{subfigure}[b]{0.4\textwidth}
	\centering
    		\includegraphics[scale=0.7]{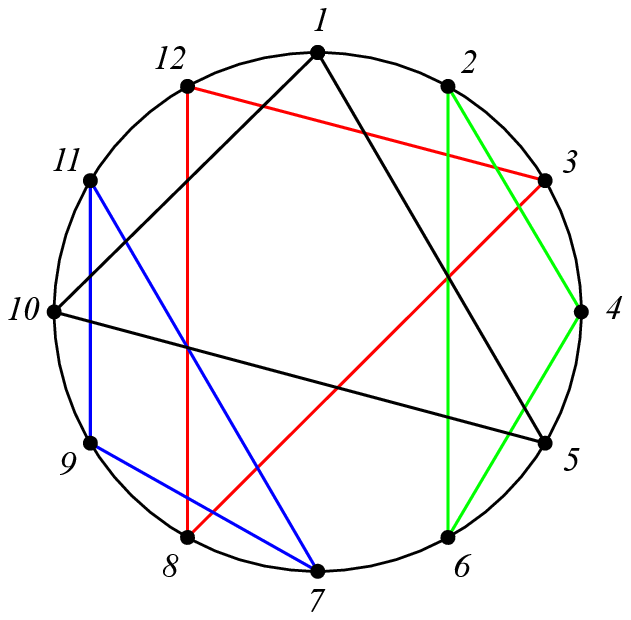}
		\caption{Generalized chord diagram}
	\end{subfigure}	
	\caption{}
\label{fig:K_3_isomorphism}
\end{figure}

In the present work we extend this approach to a more general case of an $n$-partite graph $K_{d,d,\ldots,d}$, which has $d$ vertices in each part. Any Hamiltonian cycle in such graph can be represented by a {\em generalized chord diagram} built on $n\cdot d$ vertices (Figure \ref{fig:K_3_isomorphism}). The class of such diagrams will be denoted as $B_{n}^{(d)}$. A generalized chord diagram consists of ``chords" isomorphic to graphs $K_d$ connecting $d$ points of the diagram. Similarly to the special case $d=2$, in a loopless generalized chord diagram any pair of neighboring points must belong to two different chords. 

The first part of this paper is devoted to enumerating generalized loopless chord diagrams without considering symmetries or, equivalently, to enumerating labelled Hamiltonian cycles in $K_{d,d,\ldots,d}$. The approach is based on reduction $B_{n}^{(d)}$ to so-called linear diagrams  $A_{n}^{(d)}$ \cite{Krasko_Om_chord_diagrams}. Linear diagrams have a self-contained meaning; in particular, permutations in certain classes can be depicted as linear diagrams (see, for example, \cite{Mathar}).

Depending on the notion of isomorphism used, two diagrams are said to be isomorphic if one could be obtained from the other either by a rotation or by a combination of rotations and reflections of the circle. Isomorphism classes of labelled generalized chord diagrams are said to be unlabelled generalized chord diagrams. In the second part of the paper we derive a system of recurrence relations that can be used to efficiently compute the numbers of unlabelled diagrams, and hence enumerate unlabelled Hamiltonian cycles in the graphs $K_{d,d,\ldots,d}$. We provide answers for both notions of isomorphism: for rotations only, as well as for rotations and reflections.

\section{Enumeration of labelled Hamiltonian cycles in $K_{d,d,\ldots,d}$}

As we've noted before, it will be convenient to find the numbers $b_{n}^{(d)}$ of generalized loopless chord diagrams instead of doing that for Hamiltonian cycles in the graphs $K_{d,d,\ldots,d}$ directly (Figure \ref{fig:K_3_isomorphism} (b)). 

\begin{figure}[ht]
\centering
	\centering
    	\includegraphics[scale=0.7]{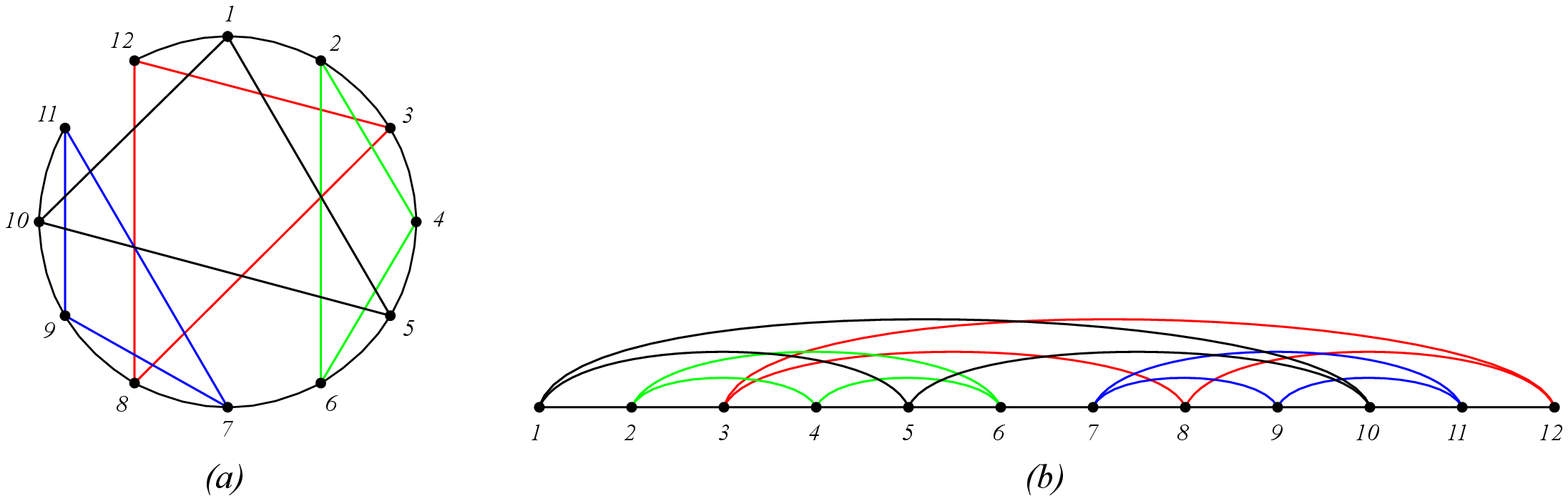}
	\caption{Cutting a triangle diagram}
\label{fig:chord_3n_cut}
\end{figure}

Each generalized chord diagram $B_n^{(d)}$ can be mapped to a unique linear diagram $A_n^{(d)}$ with $n\cdot d$ points (Figure \ref{fig:chord_3n_cut}) by cutting it (Figure \ref{fig:chord_3n_cut} (a)) along the arc that connects the points $1$ and $n\cdot d$. The result is a loopless generalized linear diagram $A_n^{(d)}$ (Figure \ref{fig:chord_3n_cut} (b)).

\begin{figure}[ht]
\centering
	\centering
    	\includegraphics[scale=0.7]{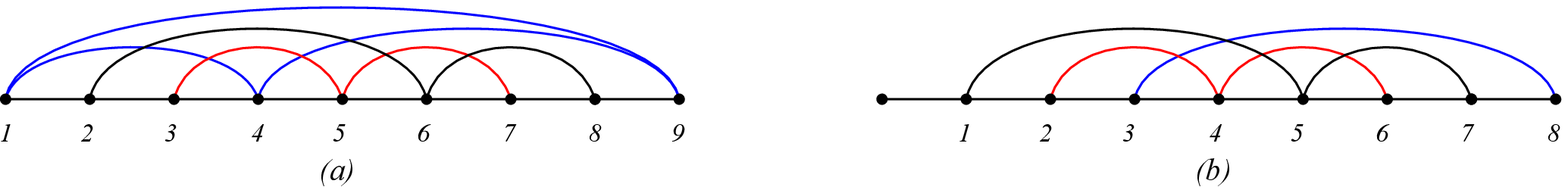}
	\caption{}
\label{fig:3_ham_b_n}
\end{figure}

Let $a_{n,k}^{(d)}$ be the number of generalized linear diagrams $A_n^{(d)}$ consisting of $d\cdot n$ points, $n$ complete subgraphs $K_d$ and having $k$ loops, $0\le k\le n(d-1)$. The numbers $b_n^{(d)}$ of generalized chord diagrams can be expressed through $a_{n,k}^{(d)}$ by the formula
\begin{equation}
\label{eq:b_d_n}
b_n^{(d)}=a_{n,0}^{(d)}-\sum\limits_{k=0}^{d-2}\BCf{d(n-1)-k-1}{d-2-k}a_{n-1,k}^{(d)}.
\end{equation}

Indeed, among all $a_{n,0}^{(d)}$ generalized linear diagrams $A_{n,0}^{(d)}$ without loops we should retain only those that have no chord connecting two end vertices. Assume that after deleting such chord $K_d$ in a diagram that contains it, we obtain a linear diagram with $k$ loops. The number of ways to obtain a generalized linear diagram $A_{n,0}^{(d)}$ from an arbitrary linear diagram $A_{n-1,k}^{(d)}$ could be counted as follows. A diagram $A_{n-1,k}^{(d)}$ has exactly $d(n-1)+1$ positions to place the vertices of the new subgraph $K_d$. We must use the first and the last of these positions, among the remaining $(d-2)$ vertices of $K_d$ we must choose some $k$ to insert them into $k$ loops, and then distribute the remaining $d-2-k$ vertices among $d(n-1)-1-k$ positions. The latter could be done in $\BCf{d(n-1)-k-1}{d-2-k}$ ways. Summing $\BCf{d(n-1)-k-1}{d-2-k}a_{n-1,k}^{(d)}$ over all possible $k$ we obtain the total number of all diagrams $A_{n,0}^{(d)}$ which have a chord connecting the first and the last point.

To find the numbers $b_n^{(d)}$ using the formula (\ref{eq:b_d_n}) we need some recurrence relations for the numbers $a_{n,k}^{(d)}$, $k=0,\ldots,d-2$. It will be easier to start with some recurrence relations for a broader range of possible values of $k$. Namely, we claim that for $k=0,\ldots,n(d-1)$ the following is true:
\begin{equation}
\label{eq:a_d_n_k}
a_{n,k}^{(d)} = \sum\limits_{t=k-d+1}^{k+d-1}{c_{n,k,t}^{(d)} \cdot a_{n-1,t}^{(d)}}, \qquad n>0,\quad 0\le k\le n(d-1),
\end{equation}
\begin{equation}
\label{eq:c_d_n_k}
c_{n,k,t}^{(d)}=\sum\limits_{i=0}^{d-1}{{d-1}\choose i}{t\choose t+i-k}{{d(n-1)-t}\choose d-2i-t+k-1},
\end{equation}
$$
a_{0,0}^{(d)}=1,\qquad \qquad a_{n,k}^{(d)}=0\qquad \text{for $n<0$, $k<0$ and $k>n(d-1)$}.
$$
The proof of relations (\ref{eq:a_d_n_k})--(\ref{eq:c_d_n_k}) is based on the procedure of removing the subgraph $K_d$ which contains the rightmost point of some generalized linear diagram $A_{n,k}^{(d)}$. Since removing $K_d$ adds or removes not more than $d-1$ loops, after removing it from some diagram  $A_{n,k}^{(d)}$ we obtain a generalized linear diagram $A_{n-1,t}^{(d)}$ with $t$ loops, $t\in[k-(d-1),k+(d-1)]$. 

\begin{figure}[ht]
\centering
	\centering
    	\includegraphics[scale=0.7]{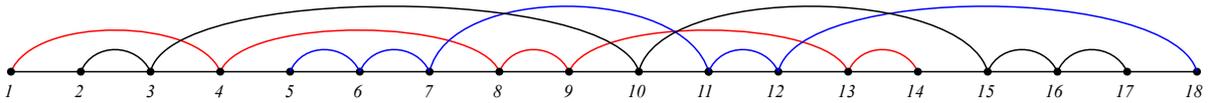}
	\caption{Generalized linear diagram $A_{3,8}^{(6)}$}
\label{fig:chord_n_k}
\end{figure}

As an example on Figure \ref{fig:chord_n_k} we show a generalized linear diagram $A_{n,k}^{(d)}$ for $n=3$, $d=6$ and $k=8$. Removing the subgraph $K_6$ that contains its rightmost point $18$ yields a generalized linear diagram $A_{2,6}^{(6)}$. 

Now assume that we are given a diagram $A_{n-1,t}^{(d)}$. The number of ways to transform it into a diagram $A_{n,k}^{(d)}$ by adding a subgraph $K_d$ which would contain its rightmost point could be counted as follows. Denote by $i$ the number of loops formed by neighboring vertices of $K_d$ after adding it to $A_{n-1,t}^{(d)}$ ($i=3$ for the subgraph $K_5$, depicted on Figure \ref{fig:chord_n_k}). For the diagram $A_{n,k}^{(d)}$ to have exactly $k$ loops after adding $K_d$, this subgraph $K_d$ must destroy $t+i-k$ existing loops ($t+i-k=1$ for the diagram $A_{2,6}^{(6)}$) by its vertices placed under these loops. Finally we have $d(n-1)-t$ remaining positions among which we can distribute $d-i-1-(t+i-k)$ remaining vertices of the subgraph $K_d$ ($6$ positions for the only vertex for the example shown on Figure \ref{fig:chord_n_k}). Counting the total number of ways to perform these combinatorial actions, we prove the formulas (\ref{eq:a_d_n_k}) -- (\ref{eq:c_d_n_k}).

Note the following special cases for the relations (\ref{eq:b_d_n})--(\ref{eq:c_d_n_k}). For $d=2$ the formula (\ref{eq:b_d_n}) becomes the formula (\ref{eq:dd_b_nk}), and the formulas (\ref{eq:a_d_n_k})--(\ref{eq:c_d_n_k}) become (\ref{eq:dd_a_nk}). For $d=3$ the expression for $b_n^{(3)}$ takes the form
$$
b_n^{(3)}=a_{n,0}-(3n-4)a_{n-1,0}^{(3)}+a_{n-1,1}^{(3)},
$$
and the system (\ref{eq:a_d_n_k})--(\ref{eq:c_d_n_k}) gets simplified to
$$
a_{n,k}=a_{n-1,k-2}+2\,[3(n-1)-(k-1)]\,a_{n-1,k-1}+\left[\BCf{3(n-1)-k}{2}+2k\right]a_{n-1,k}+
$$
$$
+(k+1)\,[3(n-1)-(k+1)]\,a_{n-1,k+1}+\BCf{k+2}{2}a_{n-1,k+2}.
$$

\section{A closed system of recurrence relations for $a_{n,k}^{(d)}$, $k=0,\ldots,d-1$} 

The recurrence relations (\ref{eq:a_d_n_k})--(\ref{eq:c_d_n_k}) in principle allow us to obtain the values of $a_{n,i}^{(d)}$, $i=0,\ldots,d-2$ that are sufficient for finding $b_n^{(d)}$. However, from the computational point of view this approach could be improved; ideally we would find a system of recurrences involving only those values of $a_{n,k}^{(d)}$ that explicitly appear in (\ref{eq:b_d_n}). For $d=2$ the approach described in \cite{Krasko_Om_chord_diagrams} was to rewrite the system (\ref{eq:a_d_n_k})--(\ref{eq:c_d_n_k}) as a system of recurrence relations, find the generating function  $w(z,t)$ for $a_{n,k}$ and then substitute $z=0$ into it. The generating function $\phi(t)=w(z,0)$ obtained as a result of substitution defines the numbers $a_{n,0}^{(2)}\equiv a_n^{(2)}$ sufficient for calculating $b_n^{(2)}$. Unfortunately this approach does not generalize well for $d>2$. One alternative approach would be to derive the corresponding system by a combinatorial argument. This approach works perfectly for $d=2$ (see \cite{Kalashnikov},\cite{Krasko_Om_chord_diagrams}), but even for $d=3$ an analogous combinatorial proof becomes quite cumbersome, and for $d>3$ the problem becomes practically intractable.

In turns out that we can actually use a combined approach: use combinatorial arguments together with the already obtained system of recurrence relations (\ref{eq:a_d_n_k})--(\ref{eq:c_d_n_k}) for the numbers $a_{n,k}^{(d)}$. With this approach we can obtain a closed system for the sequences $a_{n,k}^{(d)}$, $k=0,\ldots,d-1$, the number of which exceeds the number of terms $a_{n,k}^{(d)}$ in the formula (\ref{eq:b_d_n}) by just one. Namely, substituting  $k=0$ into the formula (\ref{eq:a_d_n_k}) we obtain the recurrence relation
\begin{equation}
\label{eq:a_n_0_d}
a_{n,0}^{(d)} = \sum\limits_{t=0}^{d-1}{c_{n,0,t}^{(d)} \cdot a_{n-1,t}^{(d)}}.
\end{equation}

\begin{figure}[ht]
\centering
	\centering
    	\includegraphics[scale=0.7]{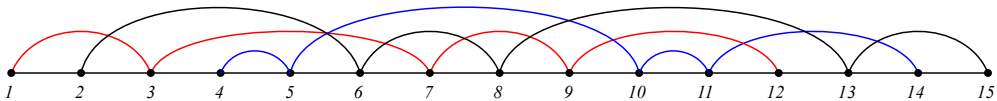}
	\caption{Generalized linear diagram $A_{3,2}^{(5)}$}
\label{fig:chord_n_k_1_a}
\end{figure}

\begin{figure}[ht]
\centering
	\centering
    	\includegraphics[scale=0.7]{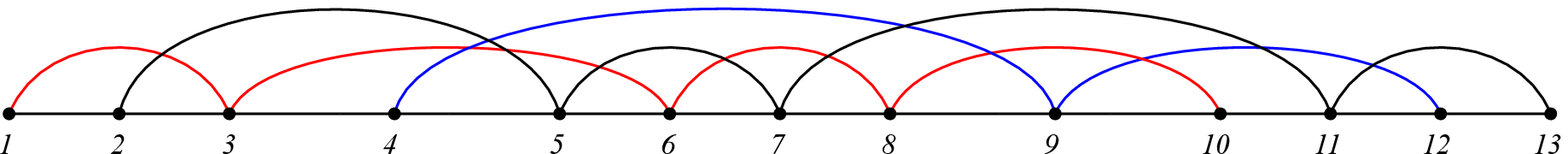}
	\caption{Reduced linear diagram}
\label{fig:chord_n_k_1_b}
\end{figure}

\begin{figure}[ht]
\centering
	\centering
    	\includegraphics[scale=0.7]{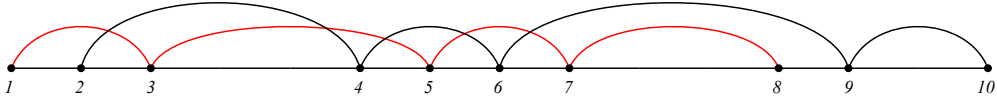}
	\caption{Generalized linear diagram $A_{2,2}^{(5)}$}
\label{fig:chord_n_k_1_c}
\end{figure}

For the values of $k$ from $1$ to $d-1$ the relations for $a_{n,k}^{(d)}$ could be obtained using combinatorial arguments. Namely, consider a generalized linear diagram $A_{n,k}^{(d)}$ that has $k$ loops distributed over $l$ subgraphs isomorphic to $K_d$, $1\leq l\leq k<d$. We will begin with the simplest case $l=1$ for which all $k$ loops are formed by a single subgraph $K_d$ (Figure \ref{fig:chord_n_k_1_a}, the corresponding subgraph $K_3$ is shown in blue). 

Contracting the loops, we transform $K_d$ into a subgraph $K_{d-k}$ in a reduced diagram which is now loopless (Figure \ref{fig:chord_n_k_1_b}). After removing this subgraph we obtain a generalized linear diagram $A_{n-1,m}^{(d)}$ with $m$ loops, $0\leq m\leq d-k$ (Figure \ref{fig:chord_n_k_1_c}, case $m=2$). Conversely, take any diagram $A_{n-1,m}^{(d)}$ and add a vertex of some new subgraph $K_{d-k}$ under $m$ of its loops. The remaining $d-k-m$ vertices of the subgraph $K_{d-k}$ should be distributed among $d(n-1)+1-m$ possible positions in $\BCf{d(n-1)+1-m}{d-k-m}$ ways. Finally, $d-k$ vertices of the subgraph $K_{d-k}$ should be transformed into a subgraph $K_d$ by replacing $k$ of its vertices with loops of the diagram. This could be done in $\BCCf{d-k}{k}=\BCf{d-1}{k}$ ways. Summing over $m$ from $0$ to $d-k$, we obtain that for $l=1$ the numbers $a_{n,k,l=1}^{(d)}$ can be expressed as
\begin{equation}
\label{eq:a_n_k_d_l=1}
a_{n,k,l=1}^{(d)}=\BCf{d-1}{k}\sum\limits_{m=0}^{d-k}\BCf{d(n-1)+1-m}{d-k-m}a_{n-1,m}^{(d)}, \qquad\qquad k=1,\ldots,d-1.
\end{equation}

\begin{figure}[ht]
\centering
	\centering
    	\includegraphics[scale=0.7]{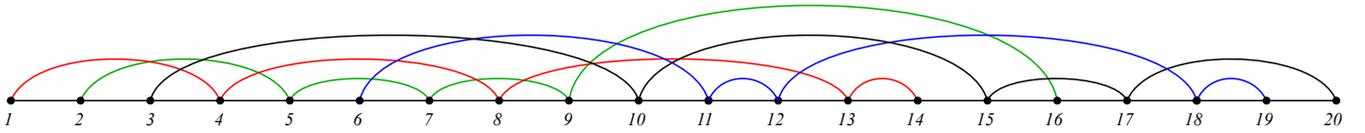}
	\caption{Generalized linear diagram $A_{4,3}^{(5)}$}
\label{fig:chord_n_k_l_a}
\end{figure}

For $l>1$ an analogous consideration becomes slightly more complicated. Indeed, let $d-r_i$, $i=1,\ldots,l$ be the number of loops in a diagram $A_{n,k}^{(d)}$ which belong to the $i$-th subgraph $K_d$, and let
\begin{equation}
\label{eq:R_set}
\sum_{i=1}^l (d-r_i)=k,\qquad 1\leq r_1\leq r_2\leq \ldots\leq r_l
\end{equation}
(see a diagram $A_{4,3}^{(5)}$ with $l=2$, $r_1=3$, $r_2=4$ on Figure \ref{fig:chord_n_k_l_a}).

\begin{figure}[ht]
\centering
	\centering
    	\includegraphics[scale=0.7]{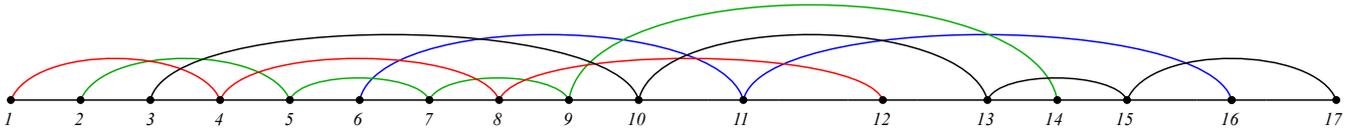}
	\caption{Reduced linear diagram}
\label{fig:chord_n_k_l_b}
\end{figure}

Contracting each such loop into a point, we obtain a reduced linear diagram with $n-l$ subgraphs $K_d$, and $l$ subgraphs $K_{r_i}$ (Figure \ref{fig:chord_n_k_l_b}). 

\begin{figure}[ht]
\centering
	\centering
    	\includegraphics[scale=0.7]{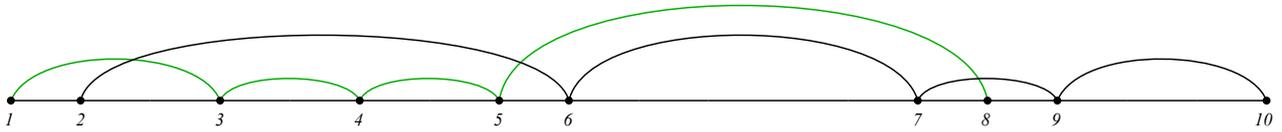}
	\caption{Generalized linear diagram $A_{2,4}^{(5)}$}
\label{fig:chord_n_k_l_c}
\end{figure}

Assume that after deleting the subgraphs $K_{r_i}$ we obtain a generalized linear diagram $A_{n-l,m}^{(d)}$, $m=0,\ldots,ld-k$ (see Figure \ref{fig:chord_n_k_l_c} corresponding to the diagram $A_{2,4}^{(5)}$). We need to determine how many diagrams $A_{n,k}^{(d)}$ with $k$ loops distributed among $l$ subgraphs $K_{r_i}$ could be obtained from this diagram $A_{n-l,m}^{(d)}$.

\begin{figure}[ht]
\centering
	\centering
    	\includegraphics[scale=0.7]{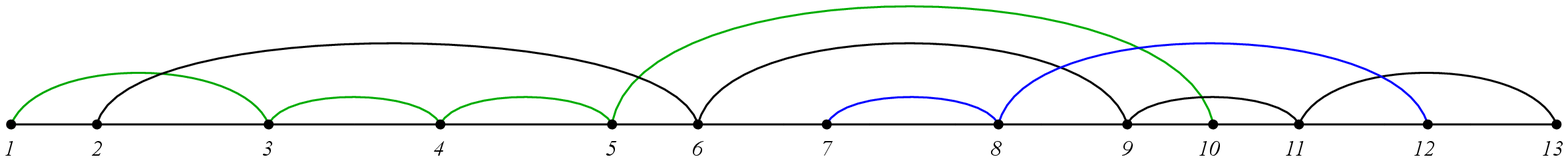}
	\caption{Generalized linear diagram $A_{2,4}^{(5)}$ with a subgraph $K_{r_1}$ added}
\label{fig:chord_n_k_l_d}
\end{figure}

Consider a diagram $A_{n-l,m_1}^{(d)}$, $m_1\equiv m$ and add a subgraph $K_{r_1}$ to it (Figure \ref{fig:chord_n_k_l_d}). Some of the existing loops may be destroyed by the vertices of $K_{r_1}$. An the same time the new diagram may have additional loops formed by neighboring vertices of the subgraph $K_{r_1}$. Denote by $s_1$ the number of loops destroyed by $K_{r_1}$, and by $j_1$ the number of loops formed by the vertices of $K_{r_1}$ ($s_1=2$, $j_1=1$ for the diagram shown on Figure \ref{fig:chord_n_k_l_d}). 

\begin{figure}[ht]
\centering
	\centering
    	\includegraphics[scale=0.7]{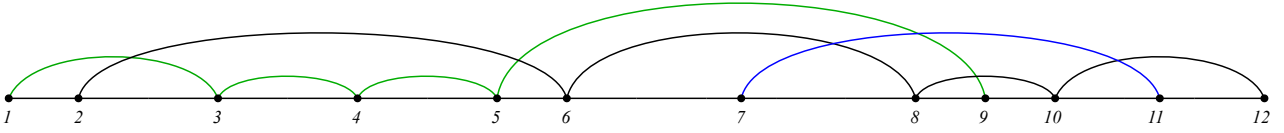}
	\caption{A generalized linear diagram $A_{2,4}^{(5)}$ with an added reduced subgraph $K_{t_1}$}
\label{fig:chord_n_k_l_e}
\end{figure}

Consider instead of $K_{r_1}$ some other subgraph $K_{t_1}$, $t_1=r_1-j_1$ obtained by contracting $j_1$ loops of the subgraph $K_{r_1}$ (Figure \ref{fig:chord_n_k_l_e}). This subgraph $K_{t_1}$ could be placed into the original diagram $A_{n-l,m_1}^{(d)}$ in such a way that $s_1$ vertices of the subgraph $K_{t_1}$ split the loops of the diagram $A_{n-l,m_1}^{(d)}$ and the remaining $t_1-s_1$ vertices are distributed among $v_1-m_1$ positions free of loops, $v_1:=d(n-l)+1$, in $\BCf{m_1}{s_1}\BCf{v_1-m_1}{t_1-s_1}$ ways. Splitting $t_1$ vertices of the subgraph $K_{t_1}$ into $r_1$ vertices such that the vertices of the obtained subgraph $K_{r_1}$ form $j_1$ additional loops could be done in $\BCCf{t_1}{j_1}=\BCCf{r_1-j_1}{j_1}=\BCf{r_1-1}{j_1}$ ways. The number
$$
\BCf{m_1}{s_1}\BCf{v_1-m_1}{r_1-j_1-s_1}\BCf{r_1-1}{j_1}
$$
of ways obtained in the first step should be multiplied by the number
$$
\BCf{m_2}{s_2}\BCf{v_2-m_2}{r_2-j_2-s_2}\BCf{r_2-1}{j_2}
$$
of ways to add a subgraph $K_{r_2}$ into $v_2:=v_1+r_1$ positions to destroy $s_2$ loops of the linear diagram with $m_2=m_1+j_1-s_1$ loops and add $j_2$ loops.

Continuing this process further, we will reach the final step where we will need to add a subgraph $K_{r_l}$ to a linear diagram. This step is special because after this addition there must be no loops in the diagram: after adding $K_{r_1},\ldots,K_{r_l}$ we must obtain a loopless reduced linear diagram (see Figure \ref{fig:chord_n_k_l_b}). Consequently, in this final step we must destroy all loops obtained on the previous step (that is, set $m_l=s_l$) and the subgraph $K_{r_l}$ should not form any loops itself (that is, $j_l=0$). 

Taking that into account, one could obtain the following final formula for the numbers $a_{n,k}^{(d)}$ for $l>1$:
\begin{equation}
\label{eq:a_n_k_d_l}
a_{n,k,l>1}^{(d)}=\sum\limits_{R}\dfrac{\alpha_R}{\beta_1!\cdot \ldots\cdot \beta_{d-1}!}
\sum\limits_{m=0}^{ld-k}p_{n,R,m}\cdot a_{n-l,m}^{(d)}.
\end{equation}
Here $R$ is an ordered multiset $\{r_1,\ldots,r_l\}$ that satisfies the conditions (\ref{eq:R_set}), the outer summation runs over all such multisets $R$ that
$$
p_{n,R,m}:=\sum\limits_{j_1=0}^{r_1-1}\ldots \sum\limits_{j_{l-1}=0}^{r_{l-1}-1}
\sum\limits_{s_1=0}^{\min\{m_1,r_1-j_1\}}\ldots \sum\limits_{s_{l-1}=0}^{\min\{m_{l-1},r_{l-1}-j_{l-1}\}}
\prod\limits_{i=1}^{l-1}\BCf{m_i}{s_i}\BCf{v_i-m_i}{r_i-j_i-s_i}\BCf{r_i-1}{j_i}\BCf{v_l-m_l}{r_l-m_l},
$$
$$
m_{i+1}:=m_i+j_i-s_i,\quad v_{i+1}:=v_i+r_i,\quad i>1;\qquad m_1:=m,\quad v_1:=d(n-l)+1,
$$
$$
\alpha_R:=\prod\limits_{i=1}^{l}\BCCf{r_i}{d-r_i}=\prod\limits_{i=1}^{l}\BCf{d-1}{r_i-1}.
$$
The multiplier $\alpha_R$ in the formula (\ref{eq:a_n_k_d_l}) describes the number of ways to transform the subgraphs $K_{r_i}$ into $K_d$. The coefficient $1/(\beta_1!\cdot \ldots\cdot\beta_{d-1}!)$ takes into account the fact that we delete the subgraphs $K_{r_i}$ not simultaneously but one after another; that is, all the cliques $K_{r_i}$ with the same number of loops are distinct. Consequently, if we have $\beta_u$ instances of a subgraph $K_u$ among all cliques $K_{r_i}$, we should divide the result by $\beta_u!$.

Finally, note that for $l>1$ the numbers $a_{n-l,m}^{(d)}$ with $m\geq d$ appear in the formula (\ref{eq:a_n_k_d_l}) for $a_{n,k}^{(d)}$. These numbers can always be eliminated using the recurrence relation  (\ref{eq:a_d_n_k}) rewritten as
\begin{equation}
\label{eq:a_n-1_k_d}
a_{n-1,k+d-1}^{(d)} = \frac{a_{n,k}^{(d)}-\sum\limits_{t=k-d+1}^{k+d-2}{c_{n,k,t}^{(d)} \cdot a_{n-1,t}^{(d)}}}{c_{n,k,k+d-1}^{(d)}}.
\end{equation}
For instance, substituting $n-1$ instead of $n$ into (\ref{eq:a_n-1_k_d}), we express the numbers $a_{n-2,d}^{(d)}$, $a_{n-2,d+1}^{(d)}$, ... through the numbers $a_{n-1,k}^{(d)}$ and $a_{n-2,m}^{(d)}$, $0\leq k,m\leq d-1$:
$$
a_{n-2,d}^{(d)} = \dfrac{a_{n-1,1}^{(d)}-\sum\limits_{t=2-d}^{d-1}{c_{n-1,1,t}^{(d)} \cdot a_{n-2,t}^{(d)}}}{c_{n-1,1,d}^{(d)}},\qquad
a_{n-2,d+1}^{(d)} = \dfrac{a_{n-1,2}^{(d)}-\sum\limits_{t=3-d}^{d}{c_{n-1,2,t}^{(d)} \cdot a_{n-2,t}^{(d)}}}{c_{n-1,2,d+1}^{(d)}},\qquad
\ldots
$$
In a similar manner we can express the numbers $a_{n-3,m}^{(d)}$, $a_{n-4,m}^{(d)}$, ... up to $a_{n-l,m}^{(d)}$.

Next we illustrate this approach using the special cases $d=2$ and $d=3$ as an example. Substituting $d=2$ into the formula (\ref{eq:a_n_0_d}), we obtain the recurrence relation of the form
$$
a_{n,0}^{(2)}=(2n-2)a_{n-1,0}^{(2)}+a_{n-1,1}^{(2)}.
$$
It can be seen that along with the numbers $a_{n,0}^{(2)}$ this equation also contains the numbers $a_{n,1}^{(2)}$ which describe linear diagrams $A_{n,1}^{(2)}$ with a single loop. For these numbers we can use the recurrence relation (\ref{eq:a_n_k_d_l=1}). Substituting the values $k=1$, $d=2$ into it, we have
$$
a_{n,1}^{(2)}=(2n-1)a_{n-1,0}^{(2)}+a_{n-1,1}^{(2)}.
$$
Expressing the numbers $a_{n,1}^{(2)}$ from these relations we obtain a second-order recurrence relation
$$
a_{n+1,0}^{(2)}=(2n+1)a_{n,0}^{(2)}+a_{n-1,0}^{(2)};\qquad\qquad a_{0,0}^{(2)}=1,\quad a_{1,0}^{(2)}=0
$$
for the number of loopless linear diagrams.

Consider a more representative example $d=3$. Substituting $d=3$ into (\ref{eq:a_n_0_d}) we have
$$
a_{n,0}^{(3)}=\BCf{3n-3}{2}a_{n-1,0}^{(3)}+(3n-4)\,a_{n-1,1}^{(3)}+a_{n-1,2}^{(3)}.
$$
The relation for $a_{n,1}^{(3)}$ as well as the recurrence relation for the numbers $a_{n,2}^{(3)}$ which corresponds to the case of both loops belonging to a single subgraph $K_3$ could be obtained from the formula (\ref{eq:a_n_k_d_l=1}): 
$$
a_{n,1}^{(3)}=2\biggl[\binom{3n-2}{2}a_{n-1,0}^{(3)}+(3n-3)\,a_{n-1,1}^{(3)}+a_{n-1,2}^{(3)}\biggr],
$$
$$
a_{n,2,l=1}^{(3)}=(3n-2)\,a_{n-1,0}^{(3)}+a_{n-1,1}^{(3)}.
$$
However, in contrast with the case $d=2$, it could happen that both loops of the diagram $A_{n,2}^{(3)}$ belong to two different subgraphs $K_3$. To count the number of such diagrams we can use the formula (\ref{eq:a_n_k_d_l}). In the current special case
$$
l=2,\qquad R=\{r_1,r_2\}=\{2,2\},\qquad \alpha_R=2\cdot 2,\qquad \beta_2=2!, \qquad j_1=j_2=1,
$$
consequently
$$
a_{n,2,l=2}^{(3)}=\frac{2\cdot 2}{2!}\sum\limits_{m=0}^{4}\sum\limits_{j=0}^{1}\sum\limits_{s=0}^{\min(m,2-j)}
\BCf{m}{s}\BCf{3(n-2)+1-m}{2-j-s}\BCf{3(n-2)+3-m-j+s}{2-m-j+s}a_{n-2,m}^{(3)}.
$$
In its turn, the numbers $a_{n-2,3}^{(3)}$ and $a_{n-2,4}^{(3)}$ are expressed through $a_{n-1,i}^{(3)}$ and $a_{n-2,j}^{(3)}$ through the relation (\ref{eq:a_n-1_k_d}):
$$
a_{n-2,3}^{(3)} = \frac{a_{n-1,1}^{(3)}-\sum\limits_{t=0}^{2}{c_{n-1,1,t}^{(3)} \cdot a_{n-2,t}^{(3)}}}{c_{n-1,3,3}^{(3)}},\qquad\qquad
a_{n-2,4}^{(3)} = \frac{a_{n-1,2}^{(3)}-\sum\limits_{t=0}^{3}{c_{n-1,2,t}^{(3)} \cdot a_{n-2,t}^{(3)}}}{c_{n-1,2,4}^{(3)}}.
$$

The obtained system of recurrence relations for $d=3$ can be simplified and rewritten as
$$
\begin{aligned}
a_{n,0}^{(3)}&=\binom {3n-3} 2 a_{n-1,0}^{(3)}+(3n-4)a_{n-1,1}^{(3)}+a_{n-1,2}^{(3)},\\[1ex]
a_{n,1}^{(3)}&=2a_{n,0}^{(3)}+2(3n-3)a_{n-1,0}^{(3)}+2a_{n-1,1}^{(3)},\\[2ex]
a_{n,2}^{(3)}&=2a_{n,0}^{(3)}+(9n-10)a_{n-1,0}^{(3)}+5a_{n-1,1}^{(3)}+2a_{n-2,0}^{(3)}.
\end{aligned}
$$

\section{Enumeration of Hamiltonian cycles in $K_{d,d,\ldots,d}$ up to rotational symmetry}

In this section we solve the problem of enumerating Hamiltonian cycles in unlabelled graphs $K_{d,d,\ldots,d}$; more precisely, we will solve an equivalent problem of enumerating unlabelled generalized chord diagrams without loops. The number $\tilde{b}_n^{(d)}$ of such diagrams can be calculated using the Burnside's lemma
\begin{equation}
\label{eq:Burnside_lemma}
\tilde{b}_n^{(d)}=\dfrac{1}{|G|}\sum\limits_{g\in G}|{\rm Fix}(g)|.
\end{equation}
Here $|\Fix(g)|$ is the number of labelled diagrams fixed by the action of an element $g$ of some group $G$ that defines the isomorphism relation between diagrams. In our case $G$ will be either the cyclic group $C_{d\cdot n}$ of diagrams' rotations or the dihedral group $D_{d\cdot n}$ of rotations and reflections. 

\begin{figure}[ht]
\centering
	\begin{subfigure}[b]{0.4\textwidth}
	\centering
    		\includegraphics[scale=1]{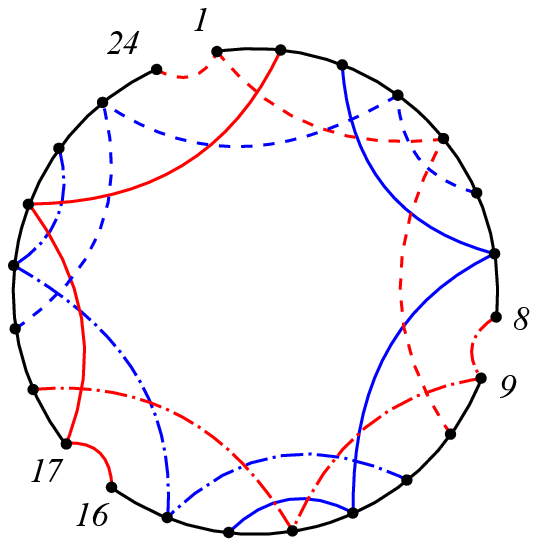}
		\caption{Diagram $A_{8,0}^{(3,4)}$}
	\end{subfigure}	
	\begin{subfigure}[b]{0.4\textwidth}
	\centering
    		\includegraphics[scale=1]{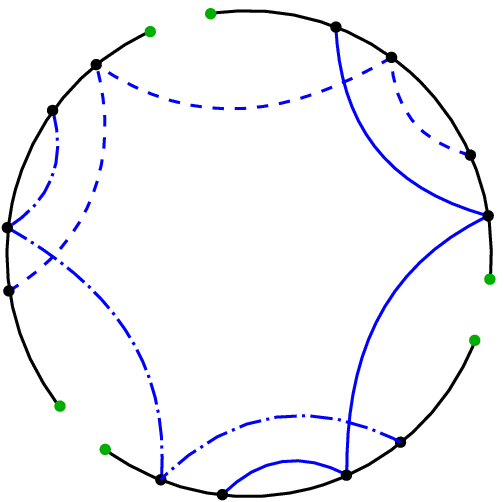}
		\caption{Diagram $A_{4,1}^{(3,4)}$}
	\end{subfigure}
	\caption{}
\label{fig:simple_diagram_K_3}
\end{figure}

\begin{figure}[ht]
\centering
	\begin{subfigure}[b]{0.4\textwidth}
	\centering
    		\includegraphics[scale=1]{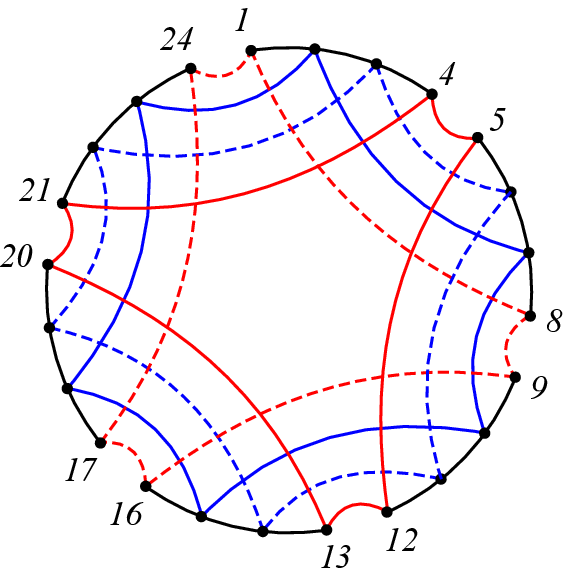}
		\caption{Diagram $A_{4,0}^{(6,6)}$}
	\end{subfigure}
	\begin{subfigure}[b]{0.4\textwidth}
	\centering
    		\includegraphics[scale=1]{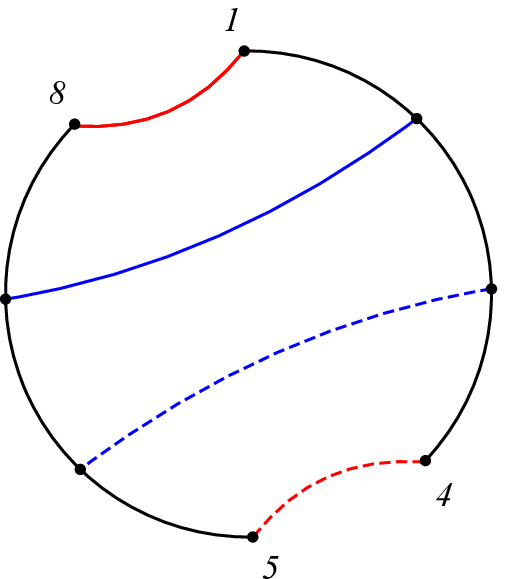}
		\caption{Reduced diagram $A_{4,0}^{(2,2)}$} 
	\end{subfigure}
	\caption{}
\label{fig:simple_diagram_K_3_b}
\end{figure}

Consider the simpler case of the cyclic group $C_{d\cdot n}$ and the action of this group on the set of generalized chord diagrams with $d\cdot n$ points and $n$ chords. Let $m$ be a divisor of $d\cdot n$, $\phi(m)$ be the Euler function of it. There are $\phi(m)$ elements of order $m$ in $C_{d\cdot n}$. Any such element fixes the same number $f(d\cdot n,m)$ of diagrams which will be called $m$-symmetric. Consequently, (\ref{eq:Burnside_lemma}) could be rewritten as
\begin{equation}
\label{eq:Burnside_lemma_Cn_d}
\tilde{b}_n^{(d)}=\dfrac{1}{d\cdot n}\sum\limits_{m\,|\, d\cdot n}\phi(m)\,f(d\cdot n,m).
\end{equation}
To calculate the values of $f(d\cdot n,m)$ it will be convenient to begin with counting so-called {\em generalized $m$-linear diagrams} (Figure \ref{fig:simple_diagram_K_3}(a) and \ref{fig:simple_diagram_K_3_b}(a)). Any such diagram with $d\cdot n$ points is obtained by cutting the circle of an $m$-symmetric generalized chord diagram into $m$ sectors between points $v$ and $v+1$, $2v$ and $2v+1$, $\ldots$, $m\cdot v$ and $1$. Each of the sectors will have $v:=d\cdot n/m$ points. By cutting between points $i$ and $i+1$ we mean that these points are no longer considered to be neighbors. Note that $1$-linear diagrams are just linear diagrams considered in the previous section. 

If we cut an $m$-symmetric loopless generalized chord diagram, we obtain a loopless generalized $m$-linear diagram. The converse is not true: if points $v$ and $v+1$, $2v$ and $2v+1$, $\ldots$, $m\cdot v$ and $1$ are connected by edges in a loopless generalized $m$-linear diagram, then gluing this diagram back into a chord diagram results in $m$ loops in it (Figure \ref{fig:simple_diagram_K_3}(a) and \ref{fig:simple_diagram_K_3_b}(a)). Denote by $A_{v,k}^{(m,d)}$ the set of generalized $m$-linear diagrams having exactly $k$ loops in each of $m$ of its sectors. The numbers $f(d\cdot n,m)$ can be expressed through the numbers $\tilde{a}_{v, k}^{(m,d)}$ of generalized $m$-linear diagrams $A_{v,k}^{(m,d)}$ by the formula
\begin{equation}
\label{eq:b_tilde_n}
f(d\cdot n,m)=\tilde{a}_{v,0}^{(m,d)}-
\sum\limits_{l\mid m}\sum\limits_{k=0}^{\tilde{d}-2}
\dfrac{l^k}{(\tilde{d}-2-k)!}\cdot\tilde{p}_{k}^{\,(\tilde{d},l)}
\cdot \tilde{a}_{v-\tilde{d},k}^{(m,d)},
\end{equation}
where $\tilde{d}=ld/m$. To prove this formula we need to show that the number of generalized loopless $m$-linear diagrams having vertices $i\cdot v$ and $i\cdot v+1$, $i=1,2,\ldots,m$ connected by edges $e_i=\{i\cdot v,i\cdot v+1\}$ is expressed by the formula
$$
\sum\limits_{l\mid m}\sum\limits_{k=0}^{\tilde{d}-2}
\dfrac{l^k}{(\tilde{d}-2-k)!}\cdot\tilde{p}_{k}^{\,(\tilde{d},l)}
\cdot \tilde{a}_{v-\tilde{d},k}^{(m,d)}.
$$ 

In the general case, $m$ edges $e_i$ may belong to $l$ different subgraphs $K_d$, $l=1,\ldots,m$ (see Figure \ref{fig:simple_diagram_K_3_b}(a)). Due to the symmetry of the diagram, $l$ must divide $m$. In addition, each clique $K_d$ covers $m/l$ loops and hence is $m/l$ symmetric. Consequently, $m/l$ must divide $d$. 

Note that it is sufficient to consider the case $l=m$ where every edge $e_i$ belongs to its own subgraph $K_d$ (see Figure \ref{fig:simple_diagram_K_3}(a)). Indeed, for any diagram with $l\neq m$ (see Figure \ref{fig:simple_diagram_K_3_b}(a)) we can select $l$ neighboring sectors ($l=2$ on Figure \ref{fig:simple_diagram_K_3_b}(a)) and build a reduced diagram with $l\cdot v$ points and $K_{\tilde{d}}$ subgraphs (see Figure \ref{fig:simple_diagram_K_3_b}(b)), in which every edge $e_i$ that connects the points $i\cdot v$ and $i\cdot v+1$ belongs to its own subgraphs $K_{\tilde{d}}$. Conversely, taking $m/l$ copies of such reduced diagrams and gluing them one after another we will obtain the diagram $A_{v,0}^{(m,d)}$. Consequently, it is sufficient to prove the formula (\ref{eq:b_tilde_n}) for the special case $l=m$ and $\tilde{d}=d$. In other words, it is sufficient to prove that the number of generalized $m$-linear diagrams in which the points $i\cdot v$ and $i\cdot v+1$, $i=1,2,\ldots,m$ are connected by edges $e_i=\{i\cdot v,i\cdot v+1\}$ each belonging to a separate subgraph $K_d$ can be calculated by the formula
\begin{equation}
\label{eq:b_tilde_n_m_d}
\sum\limits_{k=0}^{\tilde{d}-2}
\dfrac{m^k}{(d-2-k)!}\cdot\tilde{p}_{d-2-k}^{\,(d,m)}
\cdot \tilde{a}_{v-d,k}^{(l,d)}.
\end{equation}

In order to prove (\ref{eq:b_tilde_n_m_d}) consider a diagram $A_{v,0}^{(m,d)}$ and delete from it all $m$ subgraphs $K_d$ that contain edges $e_i$, $i=1,\ldots,m$. The resulting diagram has  $v-d$ vertices  and $k$ loops in each of $m$ sectors ($v-d=4$, $k=1$ on Figure \ref{fig:simple_diagram_K_3} (b)). The parameter $k$ varies from zero to $d-2$; in particular the maximal value $k=d-2$ corresponds to the case where for every point being deleted appears a new loop, unless this point is an end point of a sector. The number of ways to create a generalized $m$-linear diagram with $k$ loops is equal to $a_{v-d,k}^{(m,d)}$. We need to count the number of ways to form a chord diagram $A_{v,0}^{(m,d)}$ with edges $e_i$ connecting the vertices $i\cdot v$ and $i\cdot v+1$, $i=1,\ldots,m$ from each such diagram $A_{v-d,k}^{(m,d)}$.

Note that there are no chords in a diagram $A_{v,0}^{(m,d)}$. Hence we must insert a vertex into each of $k$ loops of the diagram $A_{v-d,k}^{(m,d)}$. That could be done in $m^k$ ways. In addition, we must insert $d-k-2$ internal (that would not lie on a sector end) vertices of the cliques $K_d$ into each sector of the diagram. We can add these vertices one after another, and that explains the divisor $(d-2-k)!$ in the denominator of the formula (\ref{eq:b_tilde_n_m_d}). In addition, we will assume that the points that fall into any of $k$ loops of the diagram are placed to the right of any $k$ points that were inserted into $k$ loops of the diagram $A_{v-d,k}^{(m,d)}$ on the first step. 

\begin{figure}[ht]
\centering
	\centering
    	\includegraphics[scale=0.7]{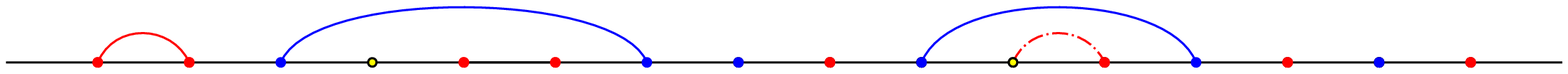}
	\caption{}
\label{fig:chord_n_k_sym}
\end{figure}

Consider a sector of a diagram $A_{v-d,k}^{(m,d)}$ (Figure \ref{fig:chord_n_k_sym}). Assume that we've already added $j$ points into it. On Figure \ref{fig:chord_n_k_sym} those $v-d$ points that belong to the original diagram $A_{v-d,k}^{(m,d)}$ are shown in blue, $k$ points inserted on the first step are shown in yellow, and $j$ additional points are shown in red. Note that the process of adding points to the diagram $A_{v-d,k}^{(m,d)}$ potentially creates some $\tilde{k}$ additional loops ($\tilde{k}=2$ on Figure \ref{fig:chord_n_k_sym}). The number $q_{\tilde{k},j+1}$ of ways to insert $j+1$ points into a sector in such a way that the number of loops $\tilde{k}$ stays the same can be expressed through the numbers $q_{\tilde{k}-1,j}$, $q_{\tilde{k},j}$ and $q_{\tilde{k}+1,j}$ by the formula
\begin{equation}
\label{eq:q_tilde_k_j}
\begin{array}{c}
q_{\tilde{k},j+1}=[2j+k-(\tilde{k}-1)]\,q_{\tilde{k}-1,j}+(\tilde{k}+1)(m-1)q_{\tilde{k}+1,j}+\\[2ex]
+[m(v-d+j+1)-(2j-\tilde{k}+k)-\tilde{k}(m-1)]q_{\tilde{k},j},\\[2ex]
q_{0,0}=1,\qquad \qquad q_{\tilde{k},j}=0\,\,\text{for $\tilde{k}<0$ or $\tilde{k}>j$}.
\end{array}
\end{equation}
Indeed, consider a diagram with $j$ red points and $\tilde{k}-1$ loops. To add a new red point into it in a way that creates one new loop, we can add it either on the left or on the right of any of $j$ red points, or to the right of any yellow point. However this counts each position that falls into one of $\tilde{k}-1$ loops twice. Consequently, the total number of positions to insert $(j+1)$-th point is equal to $2j+k-(\tilde{k}-1)$. 

To explain the multiplier $(\tilde{k}+1)(m-1)$ of $q_{\tilde{k}+1,j}$ recall that we must insert a point to destroy one of $\tilde{k}+1$ existing loops. Let such loop belong to some subgraph $K_d$. We destroy this loop in case the new point belongs to one of $m-1$ remaining complete subgraphs $K_d$. Consequently, there are $(\tilde{k}+1)(m-1)$ ways to insert this point.

Finally consider a diagram with $j$ red points and $\tilde{k}$ loops. This diagram has $v-d+j+1$ positions for inserting a new red point. More than that, this red point may belong to any of $m$ subgraphs $K_d$ that will be added. This fact explains the multiplier $m(v-d+j+1)$ in formula (\ref{eq:q_tilde_k_j}). However we should exclude those positions that lead to a new loop or to a destruction of an existing loop. A new loop can be added in $2j-\tilde{k}+k$ ways, and one of the existing loops can be destroyed in $(m-1)\tilde{k}$ ways.   

Using the recurrence relation (\ref{eq:q_tilde_k_j}) we can calculate the numbers $q_{\tilde{k},j}$ up to $j=d-2-k$. We will need the numbers $q_{0,j}$ which describe the ways to add $j$ points so that there are no loops after the addition.

The final step of building a loopless $m$-linear diagram is the addition of two end vertices to each sector. If these vertices do not belong to the same subgraph $K_d$ as their neighbors, there are no loops added and the number of ways to do that coincides with $q_{0,\tilde{v}}$, $\tilde{v}:=d-2-k$. However there exists a possibility that one or two of the end vertices creates a new loop. Hence the number $\tilde{p}_{v}$ of loopless diagrams with added end vertices can be expressed through the numbers $q_{0,j}$ using inclusion-exclusion principle:
$$
\tilde{p}_{\tilde{v}}^{\,(d,m)}=\sum\limits_{i=0}^{\tilde{v}}\sum\limits_{j=0}^{\tilde{v}-i}(-1)^{i+j}\cdot\frac{\tilde{v}!}{(\tilde{v}-i-j)!}\cdot q_{0,\tilde{v}-i-j}.
$$
In this formula the summation runs over the number $i$ of leftmost points that belong to the same clique as the point being added, and over an analogous number $j$ of rightmost points.

To use the formula (\ref{eq:b_tilde_n}) it remains to obtain some recurrence relations for the numbers $a_{v,k}^{(m,d)}$. In fact, these numbers can be counted by the following formulas:
\begin{equation}
\label{eq:a_v_k_m_d}
a_{v,k}^{(m,d)}=\sum\limits_{l\mid m}
\sum\limits_{t=k-\tilde{d}+1}^{k+\tilde{d}-1}
c_{k,t}\cdot a_{v-\tilde{d},t}^{(l,\tilde{d})},
\end{equation}
where 
\begin{equation}
\label{eq:c_k_t}
c_{k,t}= \sum\limits_{i=0}^{\tilde{d}-1}\BCf{\tilde{d}-1}{i}\cdot \BCf{t}{t+i-k}\cdot
\dfrac{l^{t+i-k}}{(\tilde{d}-2i-t+k-1)!}\cdot \hat{p}_{d-2i-t+k-1}.
\end{equation}
Note that to prove the formulas (\ref{eq:a_v_k_m_d})--(\ref{eq:c_k_t}) it is also sufficient to consider the case $l=m$, $\tilde{d}=d$: all the other cases can be reduced to it analogously as for (\ref{eq:b_tilde_n}).

In particular, to prove (\ref{eq:a_v_k_m_d}) remove $m$ cliques that contain the leftmost points of each of $m$ sectors. This yields a diagram $A_{v-d,t}^{(m,d)}$ with $v-d$ vertices and $t$ loops in each sector. The number of ways to obtain such diagram is $a_{v-d,t}^{(m,d)}$. It remains to find the numbers $c_{k,t}$ of ways to add $m$ cliques back to an arbitrary diagram $A_{v-d,t}^{(m,d)}$. 

Denote by $i$ the number of loops belonging to a clique $K_d$. Contracting such loops we obtain a clique with $d-i$ vertices. The loops can be restored in $\BCCf{d-i}{i}=\BCf{d-1}{i}$ ways. Then, similarly to the derivation of (\ref{eq:c_d_n_k}), we have to obtain $k$ loops in each sector, hence out of the existing $t$ loops we must destroy $t+i-k$. Since each loop could be destroyed by any of $m$ cliques, we should multiply the binomial coefficient $\BCf{t}{t+i-k}$ by $m^{t+i-k}$. It remains to add $v:=d-2i-t+k-1$ vertices to each sector (or to each clique). These vertices can be added one by one, and that explains the divisor $(d-2i-t+k-1)$ in the formula (\ref{eq:c_k_t}). Finally, the numbers $\hat{p}$ express the number of ways to add the vertices in such a way that the resulting diagram has no loops. The numbers $\hat{p}$ are expressed through $q_{0,j}$ by the formula
$$
\hat{p}_{\hat{v}}=\sum\limits_{i=0}^{\hat{v}}(-1)^i\cdot\frac{\hat{v}!}{(\hat{v}-i)!}\cdot q_{0,\hat{v}-i}.
$$

\section{Enumeration of Hamiltionian cycles in $K_{d,d,\ldots,d}$ up to reflections and rotations}

\begin{figure}[ht]
\centering
	\begin{subfigure}[b]{0.3\textwidth}
	\centering
    		\includegraphics[scale=1]{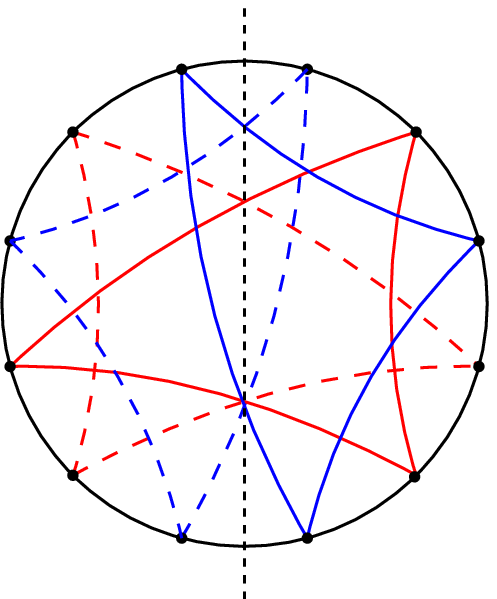}
		\caption{Diagram $H^{(0)}(4)$}
	\end{subfigure}
	\begin{subfigure}[b]{0.3\textwidth}
	\centering
    		\includegraphics[scale=1]{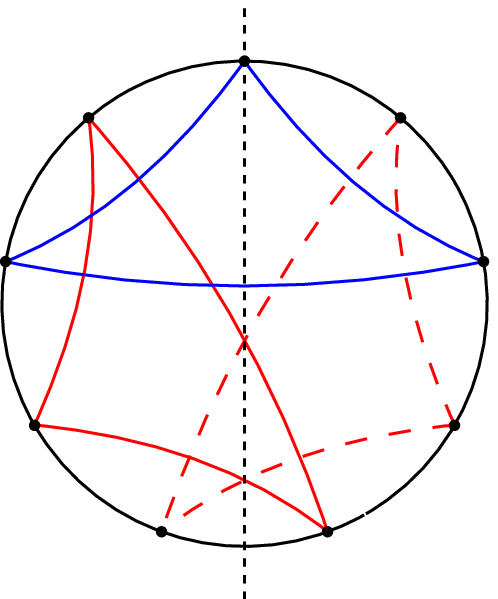}
		\caption{Diagram $H^{(1)}(3)$}
	\end{subfigure}
	\begin{subfigure}[b]{0.3\textwidth}
	\centering
    		\includegraphics[scale=1]{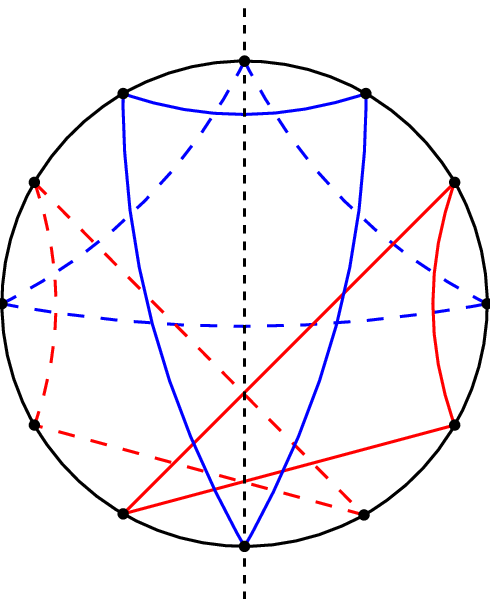}
		\caption{Diagram $H^{(2)}(4)$}
	\end{subfigure}
	\caption{}
\label{fig:refl_symmetric_odd}
\end{figure}

In this section we enumerate non-isomorphic chord diagrams under the action of the dihedral group $D_{d\cdot n}$. The Burnside lemma could be rewritten for this case as
\begin{equation}
\label{eq:Burnside_lemma_Dn}
\tilde{b}_n=\dfrac{1}{2dn}\sum\limits_{m\,|\, dn}\phi(m)\,f(dn,m)+\dfrac{h^{(0)}(n)+2h^{(1)}(n)+h^{(2)}(n)}{2},
\end{equation} 
where $h^{(i)}(n)$ denotes the number of chord diagrams symmetric under reflection about the axis passing through $i$ points of the diagram and consisting of $n$ subgraphs $K_d$ (see Figure \ref{fig:refl_symmetric_odd}).

\begin{figure}[ht]
\centering
	\begin{subfigure}[b]{0.4\textwidth}
	\centering
    		\includegraphics[scale=1]{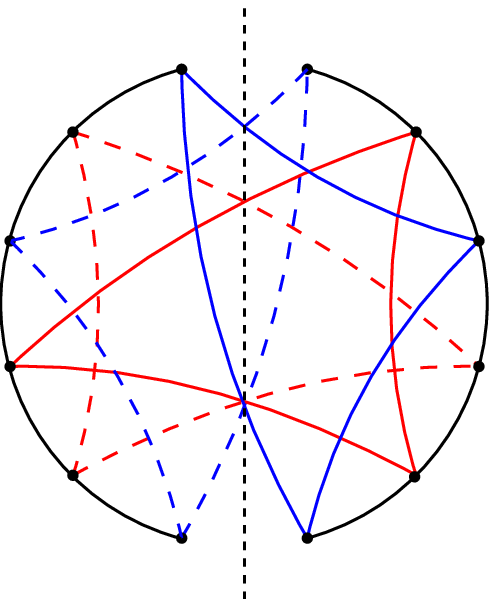}
		\caption{Diagram $H^{(0)}(4)$}
	\end{subfigure}
	\begin{subfigure}[b]{0.4\textwidth}
	\centering
    		\includegraphics[scale=1]{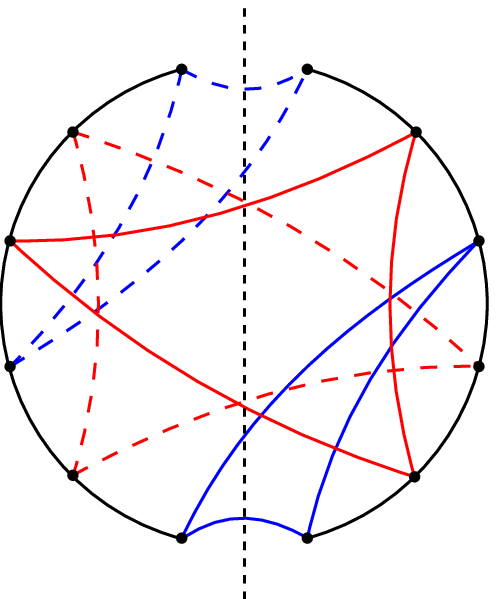}
		\caption{Diagram $A_{6,0}^{(2,3)}$} 
	\end{subfigure}
	\caption{}
\label{fig:refl_symmetric_odd_0}
\end{figure}

First consider the case of the axis of symmetry not passing through any vertices (Figure \ref{fig:refl_symmetric_odd} (a)). The number $dn$ of points must be even in this case. We can transform any such diagram into some generalized $2$-linear diagram having $dn/2$ points by cutting the circle between the points $1$ and $dn$, between the points $dn/2$ and $dn/2+1$ (Figure \ref{fig:refl_symmetric_odd_0} (a)), and then reflecting one half of the diagram along the horizontal axis (Figure \ref{fig:refl_symmetric_odd_0} (b)). However, the mapping described by this transformation may be not bijective: if a $2$-linear diagram has edges $e_1$ and $e_2$ connecting the vertices $1$ and $dn/2+1$, as well as $dn$ and $dn/2$, the reverse mapping would create two loops in the new diagram.

Note that for an odd $d$ the described mapping is actually bijective. Indeed, any edge $e_i$ must belong to some subgraph $K_d$ which must be transformed into itself by a reflection along the vertical axis, and that is impossible for odd $d$. So for odd $d$ we have
$$
h^{(0)}(n)=a_{dn/2,0}^{(2,d)}.
$$

For an even $d$ the formula is more complex: a generalized $2$-linear diagram may have both edges $e_i$ or any of them separately. Then for $h^{(0)}(n)$ we have the following formula:
$$
h^{(0)}(n)=a_{dn/2,0}^{(2,d)}-
2\sum\limits_{k=0}^{d/2-1}\alpha^{(1)}_k\cdot a_{d(n-1)/2,k}^{(2,d)}+
\sum\limits_{k=0}^{d-2}\alpha^{(2)}_k\cdot a_{d(n-2)/2,k}^{(2,d)}-
\sum\limits_{k=0}^{d/2-2}\alpha^{(3)}_k\cdot a_{d(n-2)/2,k}^{(2,d)},
$$ 
$$
\alpha^{(1)}_k:=\BCf{(n-1)d/2-1-k}{d/2-1-k},\qquad\qquad \alpha^{(3)}_k:=\BCf{(n-1)d/2-1-k}{d/2-2-k},
$$
$$
\alpha^{(2)}_k:=\sum\limits_{j=0}^{d/2-1}\sum\limits_{s=0}^{\min(k, d/2-1-j)}\BCf{d/2-1}{j}\cdot\BCf{k}{s}\cdot\BCf{d/2(n-2)-k}{d/2-1-j-s}\cdot\BCf{d/2(n-1)-1-(k+j-s)}{d/2-1-(k+j-s)}.
$$
The coefficients $\alpha^{(1)}_k$ describe the diagrams which have an edge $e$ connecting either the pair $1$ and $dn/2+1$ or the pair $dn$ and $dn/2$; the coefficients $\alpha^{(2)}_k$ describe the diagrams that have both of these pairs of vertices connected by edges $e_i$, $i=1,2$, belonging to two different subgraphs $K_d$; the coefficients $\alpha^{(3)}_k$ describe an analogous case when both $e_i$ belong to the same subgraph $K_d$.  

To find the coefficients $\alpha^{(1)}_k$ delete the subgraph $K_d$ containing the edge $e$ from a generalized $2$-linear diagram $A_{nd/2,0}^{(2,d)}$. This yields a generalized $2$-linear diagram $A_{d/2(n-1),k}^{(2,d)}$ with $k$ loops in each sector, $k=0,\ldots,d/2-2$. There are $\alpha^{(1)}_k$  ways to add a subgraph $K_d$ back to this diagram. Indeed, exactly $k$ out of $d/2-1$ points of the subgraph $K_d$ must be placed into $k$ loops of each of two sectors, one point must be placed on the end of each sector, and then the remaining $d/2-1-k$ points of the subgraph $K_d$ can be distributed among $(n-1)d/2-1$ positions in $\BCf{(n-1)d/2-1-k}{d/2-1-k}$ ways.

Consider the case of a linear diagram $A_{nd/2,0}^{(2,d)}$ having both pairs of edges connected by $e_i$. These edges could belong either to the same subgraph $K_d$ or to two different subgraphs. It would be easier to deal with the former subcase first. Removing the subgraph $K_d$ yields a diagram $A_{d/2(n-1),k}^{(2,d)}$, for which $d/2-2-k$ points of the subgraph $K_d$ must be distributed among $d/2-2-k$ positions. This could be done in $\alpha^{(3)}_k$ ways. In the latter subcase removing both subgraphs $K_d$ yields a $2$-linear diagram $A_{d/2(n-2),k}^{(2,d)}$ with $k$ loops $k=0,\ldots,d-2$. Enumerating the ways to add these two subgraphs back results in almost the same considerations as those that were performed to find the coefficients $p_{n,R,m}$. 

Namely, first we add the subgraph $K_d$ that contains the edge $e_1$ into a diagram $A_{d/2(n-2),k}^{(2,d)}$. Denote by $s$ the number of loops of the diagram that are destroyed by the vertices of $K_d$. Denote by $j$ the number of loops belonging to $K_d$ which are created in each sector. Contracting these loops we obtain a reduced subgraph with $d-2j$ vertices. One vertex of this subgraph must be placed onto the first position of a sector. Then there are $\BCf{k}{s}$ ways to place the vertices into $s$ loops and $\BCf{(n-2)d/2-k}{d/2-1-j-s}$ ways to distribute the remaining vertices over positions not covered by loops. Finally there are $\BCf{d/2-1}{j}$ ways to split $j$ vertices and transform them into the loops of $K_d$. The last binomial coefficient in $\alpha^{(2)}_k$ describes the number of ways to add the second subgraph $K_d$ so that the obtained diagram has no loops. 

The numbers $h^{(1)}_n$ describing symmetric diagrams with a single vertex lying on the axis of symmetry may be different from $0$ only for odd $n$ and $d$. Removing the subgraph $K_d$ which contains the vertex lying on the axis of symmetry and flipping one half of the diagram yields a generalized $2$-linear diagram $A_{d(n-1)/2,k}^{(2,d)}$ containing $k$ loops in each sector. Since the reverse transformation could be done in $\BCf{(n-1)d/2-1-k}{(d-1)/2-k}$ ways, the numbers $h^{(1)}(n)$ are equal to
$$
h^{(1)}(n)=\sum\limits_{k=0}^{(d-1)/2}\BCf{(n-1)d/2-1-k}{(d-1)/2-k}\cdot a_{d(n-1)/2,k}^{(2,d)}.
$$

It remains to consider the case of the axis of symmetry crossing two vertices. For an even $d$ these vertices must belong to the same clique $K_d$. Considerations analogous to those performed above for counting $h^{(1)}(n)$ yield that
$$
h^{(2)}(n)=\sum\limits_{k=0}^{d/2-1}\alpha^{(1)}_k\cdot a_{d(n-1)/2,k}^{(2,d)},\qquad\qquad \text{$d$ is even}.
$$
For an odd $d$ the numbers $h^{(2)}(n)$ may be different from $0$ only if $n$ is even. They can be calculated by the formula
$$
h^{(2)}(n)=\sum\limits_{k=0}^{d-1}\sum\limits_{j=0}^{(d-1)/2}\sum\limits_{s=0}^{\min(k, (d-1)/2-j)}
\BCf{(d-1)/2}{j}\BCf{k}{s}\BCf{d/2(n-2)-k}{(d-1)/2-j-s}\times
$$
$$
\times\BCf{d/2(n-1)-1-(k+j-s)}{(d-1)/2-(k+j-s)}\cdot a_{d(n-2)/2,k}^{(2,d)}.
$$

\section*{Conclusion}
In this paper labelled and unlabelled generalized loopless chord and linear diagrams were enumerated. Labelled chord diagrams can be thought of as directed Hamiltonian cycles with a distinguished starting point in unlabelled graphs $K_{d,d,\ldots,d}$. Considering chord diagrams up to rotations removes the starting point but keeps the direction chosen. Finally, chord diagrams considered up to rotations and reflections is nothing but undirected unlabelled Hamiltonian cycles in unlabelled graphs $K_{d,d,\ldots,d}$. The final numbers for various classes of linear and chord diagrams as well as Hamiltonian cycles considered above can be found in Tables \ref{table:loopless1}--\ref{table:loopless5}.

\begin{table}[h!]
\scriptsize
\centering
\begin{tabular}{c|cccc}
\midrule
$ n $  &\phantom{0}Linear\phantom{0}&\phantom{0}Chord labelled\phantom{0}&\phantom{0}Up to rotations\phantom{0}&\phantom{0}Up to all symmetries\phantom{0}\\ 
\midrule
1 & 0 & 0 & 0 & 0 \\
2 & 1 & 1 & 1 & 1 \\
3 & 5 & 4 & 2 & 2 \\
4 & 36 & 31 & 7 & 7 \\
5 & 329 & 293 & 36 & 29 \\
6 & 3655 & 3326 & 300 & 196 \\
7 & 47844 & 44189 & 3218 & 1788 \\
8 & 721315 & 673471 & 42335 & 21994 \\
9 & 12310199 & 11588884 & 644808 & 326115 \\
10 & 234615096 & 222304897 & 11119515 & 5578431 \\
11 & 4939227215 & 4704612119 & 213865382 & 107026037 \\
12 & 113836841041 & 108897613826 & 4537496680 & 2269254616 \\
13 & 2850860253240 & 2737023412199 & 105270612952 & 52638064494 \\
14 & 77087063678521 & 74236203425281 & 2651295555949 & 1325663757897 \\
15 & 2238375706930349 & 2161288643251828 & 72042968876506 & 36021577975918 \\
16 & 69466733978519340 & 67228358271588991 & 2100886276796969 & 1050443713185782 \\
17 & 2294640596998068569 & 2225173863019549229 & 65446290562491916 & 32723148860301935 \\
18 & 80381887628910919255 & 78087247031912850686 & 2169090198219290966 & 1084545122297249077 \\
19 & 2976424482866702081004 & 2896042595237791161749 & 76211647261082309466 & 38105823782987999742 \\
20 & 116160936719430292078411 & 113184512236563589997407 & 2829612806029873399561 & 1414806404051118314077 \\
\midrule\end{tabular}
\caption{Loopless diagrams by number of $K_2$}
\label{table:loopless1}
\end{table}

\begin{table}[h!]
\scriptsize
\centering
\begin{tabular}{c|cccc}
\midrule
$ n $  &\phantom{0}Linear\phantom{0}&\phantom{0}Chord labelled\phantom{0}&\phantom{0}Up to rotations\phantom{0}&\phantom{0}Up to all symmetries\phantom{0}\\ 
\midrule
1 & 0 & 0 & 0 & 0 \\
2 & 1 & 1 & 1 & 1 \\
3 & 29 & 22 & 4 & 4 \\
4 & 1721 & 1415 & 126 & 83 \\
5 & 163386 & 140343 & 9367 & 4848 \\
6 & 22831355 & 20167651 & 1120780 & 562713 \\
7 & 4420321081 & 3980871156 & 189565588 & 94810999 \\
8 & 1133879136649 & 1035707510307 & 43154533233 & 21577786374 \\
9 & 372419001449076 & 343866839138005 & 12735808866899 & 6367912802891 \\
10 & 152466248712342181 & 141979144588872613 & 4732638168795171 & 2366319275431001 \\
11 & 76134462292157828285 & 71386289535825383386 & 2163220895025390670 & 1081610451348718567 \\
12 & 45552714996556390334921 & 42954342000612934599071 & 1193176166690983987122 & 596588083450068950934 \\
13 & 32173493282909179882613934 & 30482693813120122213093587 & 781607533669746761791541 & 390803766837390136477505 \\
\midrule\end{tabular}
\caption{Loopless diagrams by number of $K_3$}
\label{table:loopless2}
\end{table}

\begin{table}[h!]
\scriptsize
\centering
\begin{tabular}{c|cccc}
\midrule
$ n $  &\phantom{0}Linear\phantom{0}&\phantom{0}Chord labelled\phantom{0}&\phantom{0}Up to rotations\phantom{0}&\phantom{0}Up to all symmetries\phantom{0}\\ 
\midrule
1 & 0 & 0 & 0 & 0 \\
2 & 1 & 1 & 1 & 1 \\
3 & 182 & 134 & 15 & 13 \\
4 & 94376 & 75843 & 4790 & 2576 \\
5 & 98371884 & 83002866 & 4151415 & 2081393 \\
6 & 182502973885 & 158861646466 & 6619291247 & 3309962320 \\
7 & 551248360550999 & 490294453324924 & 17510518983528 & 8755277273334 \\
8 & 2536823683737613858 & 2292204611710892971 & 71631394311300461 & 35815698613833466 \\
9 & 16904301142107043464659 & 15459367618357013402267 & 429426878302882412435 & 214713439275724149414 \\
10 & 156690501089429126239232946 & 144663877588996810362218074 & 3616596939726424941979785 & 1808298469877117320495867 \\
\midrule\end{tabular}
\caption{Loopless diagrams by number of $K_4$}
\label{table:loopless3}
\end{table}

\begin{table}[h!]
\scriptsize
\centering
\begin{tabular}{c|cccc}
\midrule
$ n $  &\phantom{0}Linear\phantom{0}&\phantom{0}Chord labelled\phantom{0}&\phantom{0}Up to rotations\phantom{0}&\phantom{0}Up to all symmetries\phantom{0}\\ 
\midrule
1 & 0 & 0 & 0 & 0 \\
2 & 1 & 1 & 1 & 1 \\
3 & 1198 & 866 & 60 & 42 \\
4 & 5609649 & 4446741 & 222477 & 112418 \\
5 & 66218360625 & 55279816356 & 2211192688 & 1105696796 \\
6 & 1681287695542855 & 1450728060971387 & 48357603758012 & 24178822553773 \\
7 & 81644850343968535401 & 72078730629785795963 & 2059392303708166507 & 1029696155560021174 \\
8 & 6945222145021508480249929 & 6235048155225093080061949 & 155876203880714141444480 & 77938101941693076258854 \\
\midrule\end{tabular}
\caption{Loopless diagrams by number of $K_5$}
\label{table:loopless4}
\end{table}

\begin{table}[h!]
\scriptsize
\centering
\begin{tabular}{c|cccc}
\midrule
$ n $  &\phantom{0}Linear\phantom{0}&\phantom{0}Chord labelled\phantom{0}&\phantom{0}Up to rotations\phantom{0}&\phantom{0}Up to all symmetries\phantom{0}\\ 
\midrule
1 & 0 & 0 & 0 & 0 \\
2 & 1 & 1 & 1 & 1 \\
3 & 8142 & 5812 & 335 & 203 \\
4 & 351574834 & 276154969 & 11508322 & 5765385 \\
5 & 47940557125969 & 39738077935264 & 1324603148183 & 662305416760 \\
6 & 16985819072511102549 & 14571371516350429940 & 404760320241653655 & 202380163158922023 \\
7 & 13519747358522016160671387 & 11876790400066163254723167 & 282780723811372935744420 & 141390361908351519807928 \\
\midrule\end{tabular}
\caption{Loopless diagrams by number of $K_6$}
\label{table:loopless5}
\end{table}

\section*{Acknowledgments}

The research was supported by the Russian Foundation for Basic Research (grant 17-01-00212).

\end{document}